\documentclass[12pt]{amsart}
\usepackage[francais,english]{babel}

\usepackage[]{amsfonts}
\usepackage{amssymb}
\usepackage{amsmath}

\def\couleur(#1 #2 #3)
	{
\special{ps::[end]}}

\def\underset#1#2{\mathrel{\mathop{\kern0pt #2}\limits_{#1}}}

\def\overset#1#2{\mathrel{\mathop{\kern0pt #2}\limits^{#1}}}

\def\bx#1{\setbox1=\hbox{\kern3pt{#1}\kern3pt}			
 \dimen1=\ht1 \advance\dimen1 by 3pt \dimen2=\dp1 \advance\dimen2 by 3pt
 \setbox1=\hbox{\vrule height\dimen1 depth\dimen2\box1\vrule}%
 \setbox1=\vbox{\hrule\box1\hrule}%
 \advance\dimen1 by .4pt \ht1=\dimen1
 \advance\dimen2 by .4pt \dp1=\dimen2 \box1\relax}

\def\wbb#1{\kern#1em}
\def\vci{\vrule  width.02em height1.47ex depth-.0ex}		
\def\11{{\rm\wbb{.2}\vci\wbb{-.37}1}}

\newcommand{\Supp}{\operatorname{Supp}} 

\newtheorem{Theorem}{Theorem}[section]
\newtheorem{Lemma}[Theorem]{Lemma}

\begin{document}

\title{Serre duality and H\"ormander's solution of the  $\bar \partial $ -equation.}

\author{Eric Amar}
\maketitle
 \ \par 
\ \par 
\renewcommand{\abstractname}{Abstract}

\begin{abstract}
We use duality in the manner of Serre to generalize a theorem
 of Hedenmalm on solution of the  $\bar \partial $  equation
 with inverse of the weight in H\"ormander  $\displaystyle L^{2}$
  estimates.\ \par 
\end{abstract}

\section{Introduction}
\quad \quad  	Let  $\displaystyle \varphi $  be a  ${\mathcal{C}}^{2}$  strictly
 sub harmonic function in the complex plane  ${\mathbb{C}},$
  i.e. with  $\Delta $  the Laplacian,  $\Delta \varphi >0$ 
 in  ${\mathbb{C}}.$  Let  $A^{2}({\mathbb{C}},e^{-2\varphi })$
  be the set of all holomorphic functions  $g$  in  ${\mathbb{C}}$
  such that\ \par 
\quad \quad \quad \quad \quad 	 $\ {\left\Vert{g}\right\Vert}_{L^{2}({\mathbb{C}},e^{-2\varphi
 })}:=\int_{{\mathbb{C}}}{\left\vert{g}\right\vert ^{2}e^{-2\varphi
 }dA}<\infty ,$ \ \par 
with  $\displaystyle dA$  the Lebesgue measure in  ${\mathbb{C}}.$
  Suppose that  $f\in L^{2}({\mathbb{C}},e^{2\varphi })$  verifies\ \par 
\quad \quad \quad \quad  	\begin{equation}  \forall g\in A^{2}({\mathbb{C}},e^{-2\varphi
 }),\ \int_{{\mathbb{C}}}{fgdA}=0\label{LpPoidsStein12}\end{equation}\ \par 
then 	in a recent paper H. Hedenmalm~\cite{Hedenmalm13} proved\ \par 
\begin{Theorem}
 Suppose that  $\displaystyle f\in L^{2}({\mathbb{C}},e^{2\varphi
 })$  verifies condition~(\ref{LpPoidsStein12}) then there exists
 a solution to the  $\bar \partial $ -equation  $\bar \partial u=f$  with\par 
\quad \quad \quad \quad \quad 	 $\ \int_{{\mathbb{C}}}{\left\vert{u}\right\vert ^{2}e^{2\varphi
 }\Delta \varphi dA}\leq \frac{1}{2}\int_{{\mathbb{C}}}{\left\vert{f}\right\vert
 ^{2}e^{2\varphi }dA}.$ \par 
\end{Theorem}
He adds in remark 1.3. that this theorem should generalize to
 the setting of several complex variables. The aim of this note
 is to show that he was right and in fact we have to make slight
 modifications in our paper~\cite{AnGrauLr12}, inspired by Serre's
 duality theorem~\cite{Serre55} to get it. The paper~\cite{AnGrauLr12}
 contains more material and here we extract just the part we
 need to make this work self contained.\ \par 
\ \par 
\quad \quad  	Let  $\varphi $  be a strictly plurisubharmonic function of
 class  ${\mathcal{C}}^{2}$  in the Stein manifold  $\displaystyle
 \Omega .$  Let  $\displaystyle c_{\varphi }(z)$  be the smallest
 eigenvalue of  $\partial \bar \partial \varphi (z),$  then 
 $\displaystyle \forall z\in \Omega ,\ c_{\varphi }(z)>0.$  We
 denote by  $\displaystyle L^{2}_{p,q}(\Omega ,e^{\varphi })$
  the set of  $\displaystyle (p,q)$  currents  $\omega $  whose
 coefficients are in the space  $\displaystyle L^{2}(\Omega ,e^{\varphi
 }),$  i.e. there is a constant  $\displaystyle C>0$  such that
 in a coordinates patch  $\displaystyle (U,\psi )$ , we have,
 with  $\displaystyle dm$  the Lebesgue measure in  ${\mathbb{C}}^{n},$ \ \par 
\quad \quad \quad \quad   $\psi ^{*}\omega =\sum_{I,J}{\omega _{I,J}dz^{I}d\bar z^{J}},\
 \int_{\psi (U)}{\left\vert{\omega _{I,J}}\right\vert ^{2}e^{\varphi
 \circ \psi ^{-1}dm}}\leq C.$ 	 $\displaystyle $ \ \par 
We denote by  $\displaystyle L^{2,c}_{p,q}(\Omega ,e^{\varphi
 })$  the currents in  $\displaystyle L^{2}_{p,q}(\Omega ,e^{\varphi
 })$  with {\sl compact support} in  $\displaystyle \Omega $
  and   ${\mathcal{H}}_{p}(\Omega )$  the set of all  $(p,\ 0),\
 \bar \partial $  closed forms in  $\Omega .$  If  $p=0,\ {\mathcal{H}}_{0}(\Omega
 )={\mathcal{H}}(\Omega )$   is the set of holomorphic functions
 in  $\Omega .$ \ \par 
\quad \quad  	We shall prove\ \par 
\begin{Theorem}
~\label{LpPoidsStein11}Let  $\displaystyle \Omega $  be a pseudo
 convex domain in  ${\mathbb{C}}^{n}$ ; if  $\omega \in L^{2,c}_{p,q}(\Omega
 ,\ e^{\varphi })$  with  $\displaystyle \bar \partial \omega
 =0$  if  $\displaystyle q<n$  and  $\displaystyle \omega \in
 L^{2}_{p,q}(\Omega ,e^{\varphi })$  with  $\displaystyle \omega
 \perp {\mathcal{H}}^{p}(\Omega )$  if  $\displaystyle q=n,$
  then there is  $\displaystyle u\in L^{2}_{p,q-1}(\Omega ,c_{\varphi
 }e^{\varphi })$  such that  $\displaystyle \bar \partial u=\omega ,$  and\par 
\quad \quad \quad \quad \quad \quad \quad 	 $\displaystyle \ {\left\Vert{u}\right\Vert}_{L^{2}(\Omega ,c_{\varphi
 }e^{\varphi })}\leq C{\left\Vert{\omega }\right\Vert}_{L^{2}(\Omega
 ,e^{\varphi })}.$ \par 
\end{Theorem}
\quad \quad  	In the case  $\displaystyle \Omega $  is a Stein manifold,
 the result is more restrictive :\ \par 
\begin{Theorem}
~\label{LpPoidsStein13}Let  $\displaystyle \Omega $  be a Stein
 manifold ; there is a convex increasing function  $\chi $  such
 that, with  $\displaystyle \psi :=\chi (\varphi ),$  if  $\omega
 \in L^{2,c}_{p,q}(\Omega ,\ e^{\psi })$  with  $\displaystyle
 \bar \partial \omega =0$  if  $\displaystyle q<n$  and  $\displaystyle
 \omega \in L^{2}_{p,q}(\Omega ,e^{\psi })$  with  $\displaystyle
 \omega \perp {\mathcal{H}}^{p}(\Omega )$  if  $\displaystyle
 q=n,$  then there is  $\displaystyle u\in L^{2}_{p,q-1}(\Omega
 ,c_{\psi }e^{\psi })$  such that  $\displaystyle \bar \partial
 u=\omega ,$  and\par 
\quad \quad \quad \quad \quad \quad \quad 	 $\displaystyle \ {\left\Vert{u}\right\Vert}_{L^{2}(\Omega ,c_{\psi
 }e^{\psi })}\leq C{\left\Vert{\omega }\right\Vert}_{L^{2}(\Omega
 ,e^{\psi })}.$ \par 
\end{Theorem}
\quad \quad  	Clearly theorem~\ref{LpPoidsStein11} generalizes Hedenmalm's
 theorem, because in one variable, we have  $\displaystyle q=n=1$
  and no compactness assumption is required.\ \par 

\section{The proof.}
\quad \quad  	Set a weight  $\eta :=e^{-\varphi }$  and  $\displaystyle \mu
 :=c_{\varphi }^{-1}\eta $  such that if  $\displaystyle \alpha
 \in L_{(p,q)}^{2}(\Omega ,\mu )$  is such that  $\bar \partial
 \alpha =0$  in  $\displaystyle \Omega $  then there is  $\displaystyle
 (p,q-1)$  current  $\displaystyle \varphi \in L^{2}(\Omega ,\eta
 )$  with  $\bar \partial \varphi =\alpha .$  Moreover we want
 that if  $\beta \in L_{(p,q)}^{2}(\Omega ,\eta ),\ \bar \partial
 \beta =0,$  there is a  $\displaystyle \gamma \in L_{(p,q-1)}^{2}(\Omega
 ,loc)$  such that  $\bar \partial \gamma =\beta .$ \ \par 
\quad \quad  	Now we suppose that, if  $\displaystyle q<n$  then  $\omega
 $  has a compact support,  $\omega \in L_{(p,q)}^{2,c}(\Omega
 ,\eta ^{-1})$  and  $\bar \partial \omega =0$  and if  $\displaystyle
 q=n,$  then  $\displaystyle \omega \perp {\mathcal{H}}_{n-p}(\Omega
 )$  and  $\displaystyle \omega \in L_{(p,q)}^{2}(\Omega ,\eta
 ^{-1}),$  but  $\omega $  needs not have compact support.\ \par 
 We copy lemma 3.5 from~\cite{AnGrauLr12}.\ \par 
\begin{Lemma}
 ~\label{GraAnd3}For  $\displaystyle \Omega ,\ \omega ,\ \eta
 $  as above, the function  ${\mathcal{L}}={\mathcal{L}}_{\omega
 },$  defined on  $(n-p,n-q+1)$  form  $\alpha \in L^{2}(\Omega
 ,\beta ),\ \bar \partial $  closed in  $\Omega ,$  as follows :\par 
\quad \quad \quad \quad \quad 	 ${\mathcal{L}}(\alpha ):=(-1)^{p+q-1}{\left\langle{\omega ,\varphi
 }\right\rangle},$  where  $\varphi \in L^{2}(\Omega ,\eta )$
  is such that  $\bar \partial \varphi =\alpha $  in  $\Omega $ \par 
is well defined and linear.\par 
\end{Lemma}
\quad \quad  	Proof.\ \par 
We have  $\alpha \in L^{2}(\Omega ,\mu ),\bar \partial \alpha
 =0,$  then such a  $\varphi \in L^{2}(\Omega ,\eta )$  exists
 by hypothesis and the pairing  $\displaystyle \ {\left\langle{\omega
 ,\varphi }\right\rangle}$  is well defined because  $\omega
 \in L_{(p,q)}^{2}(\Omega ,\eta ^{-1}).$ \ \par 
Let us see that  ${\mathcal{L}}$  is well defined.\ \par 
\quad \quad  	Suppose first that  $q<n.$ \ \par 
In order for  ${\mathcal{L}}$  to be well defined we need\ \par 
\quad \quad \quad   $\forall \varphi ,\psi \in L^{2}_{(n-p,n-q)}(\Omega ,\eta ),\
 \bar \partial \varphi =\bar \partial \psi \Rightarrow {\left\langle{\omega
 ,\varphi }\right\rangle}={\left\langle{\omega ,\psi }\right\rangle}.$ \ \par 
Then we have  $\bar \partial (\varphi -\psi )=0$  hence we can
 solve  $\bar \partial $  in  $L^{2}(\Omega ,loc)$  :\ \par 
\quad \quad \quad \quad \quad 	 $\exists \gamma \in L^{2}_{(n-p,n-q-1)}(\Omega ,loc)::\bar
 \partial \gamma =(\varphi -\psi ).$ \ \par 
So  $\ {\left\langle{\omega ,\varphi -\psi }\right\rangle}={\left\langle{\omega
 ,\bar \partial \gamma }\right\rangle}=(-1)^{p+q-1}{\left\langle{\bar
 \partial \omega ,\gamma }\right\rangle}=0$  because  $\omega
 $  is compactly supported in  $\Omega .$ \ \par 
Hence  ${\mathcal{L}}$  is well defined in that case.\ \par 
\quad \quad  	Suppose now that  $q=n.$ \ \par 
Of course  $\bar \partial \omega =0$  and we have that  $\varphi
 ,\ \psi $  are  $(n-p,\ 0)$  currents hence  $\bar \partial
 (\varphi -\psi )=0$  means that  $h:=\varphi -\psi $  is a 
 $\bar \partial $  closed  $(n-p,\ 0)$  current hence  $h\in
 {\mathcal{H}}_{n-p}(\Omega ).$  Hence the hypothesis,  $\displaystyle
 \omega \perp {\mathcal{H}}_{n-p}(\Omega )$  gives  $\ {\left\langle{\omega
 ,h}\right\rangle}=0,$  and  ${\mathcal{L}}$  is also well defined
 in that case and no compactness on the support of  $\omega $
  is required here.\ \par 
\quad \quad  	It remains to see that  ${\mathcal{L}}$  is linear, so let
  $\alpha =\alpha _{1}+\alpha _{2},$  with  $\displaystyle \alpha
 _{j}\in L^{2}(\Omega ,\mu ),\ \bar \partial \alpha _{j}=0,\
 j=1,2\ ;$  we have  $\alpha =\bar \partial \varphi ,\ \alpha
 _{1}=\bar \partial \varphi _{1}$  and  $\alpha _{2}=\bar \partial
 \varphi _{2},$  with  $\varphi ,\ \varphi _{1},\ \varphi _{2}$
  in  $\displaystyle L^{2}(\Omega ,\eta )$  so, because  $\bar
 \partial (\varphi -\varphi _{1}-\varphi _{2})=0,$  we have\ \par 
 $\bullet $ {\it  }for  $\displaystyle q<n,$ \ \par 
\quad \quad \quad 	 $\varphi =\varphi _{1}+\varphi _{2}+\bar \partial \psi ,$ 
 with  $\psi $  in  $L^{2}(\Omega ,loc),$ \ \par 
so\ \par 
\quad \quad \quad 	 ${\mathcal{L}}(\alpha )=(-1)^{p+q-1}{\left\langle{\omega ,\varphi
 }\right\rangle}=(-1)^{p+q-1}{\left\langle{\omega ,\varphi _{1}+\varphi
 _{2}+\bar \partial \psi }\right\rangle}=\ {\mathcal{L}}(\alpha
 _{1})+{\mathcal{L}}(\alpha _{2})+(-1)^{p+q-1}{\left\langle{\omega
 ,\bar \partial \psi }\right\rangle},$ \ \par 
but  $\ (-1)^{p+q-1}{\left\langle{\omega ,\bar \partial \psi
 }\right\rangle}={\left\langle{\bar \partial \omega ,\psi }\right\rangle}=0,$
  because  $\displaystyle \Supp \omega \Subset \Omega $  implies
 there is no boundary term  so  $\displaystyle {\mathcal{L}}(\alpha
 )={\mathcal{L}}(\alpha _{1})+{\mathcal{L}}(\alpha _{2}).$ \ \par 
 $\bullet $  for  $\displaystyle q=n,$ \ \par 
\quad \quad 	 $\displaystyle \varphi =\varphi _{1}+\varphi _{2}+h,$  with
  $\displaystyle h\in {\mathcal{H}}_{n-p}(\Omega )$  hence\ \par 
\quad \quad \quad \quad \quad 	 $\displaystyle {\mathcal{L}}(\alpha )=(-1)^{p+q-1}{\left\langle{\omega
 ,\varphi }\right\rangle}=(-1)^{p+q-1}{\left\langle{\omega ,\varphi
 _{1}+\varphi _{2}+h}\right\rangle}=\ {\mathcal{L}}(\alpha _{1})+{\mathcal{L}}(\alpha
 _{2})+(-1)^{p+q-1}{\left\langle{\omega ,h}\right\rangle},$ \ \par 
so, because  $\displaystyle \ {\left\langle{\omega ,h}\right\rangle}=0,$
  we still have  $\displaystyle {\mathcal{L}}(\alpha )={\mathcal{L}}(\alpha
 _{1})+{\mathcal{L}}(\alpha _{2})$  without compactness assumption
 on the support of  $\omega .$ \ \par 
The same for  $\alpha =\lambda \alpha _{1}$  and the linearity.
  $\hfill\blacksquare $ \ \par 
\ \par 
\begin{Lemma}
~\label{GraAnd4}Still with the same hypotheses as above there
 is a  $\displaystyle (p,q-1)$  current  $u$  such that\par 
\quad \quad 	 $\displaystyle \forall \alpha \in L_{(n-p,n-q+1)}^{2}(\Omega
 ,\mu ),\ {\left\langle{u,\alpha }\right\rangle}={\mathcal{L}}(\alpha
 )=(-1)^{p+q-1}{\left\langle{\omega ,\varphi }\right\rangle},$ \par 
and\par 
\quad \quad \quad \quad \quad 	 $\displaystyle \sup _{\alpha \in L^{2}(\Omega ,\mu ),\ {\left\Vert{\alpha
 }\right\Vert}_{L^{2}(\Omega ,\mu )}\leq 1}\ \left\vert{{\left\langle{u,\alpha
 }\right\rangle}}\right\vert \leq C{\left\Vert{\omega }\right\Vert}_{L^{2}(\Omega
 ,\eta )}.$ \par 
\end{Lemma}
\quad \quad  	Proof.\ \par 
By lemma~\ref{GraAnd3} we have that  ${\mathcal{L}}$  is a linear
 form on  $(n-p,n-q+1)$  forms  $\alpha \in L^{}(\Omega ,\mu
 ),\ \bar \partial $  closed in  $\Omega .$ \ \par 
	We have\ \par 
\quad \quad \quad \quad \quad 	 $\exists \varphi \in L_{(n-p,n-q)}^{2}(\Omega ,\eta )::\bar
 \partial \varphi =\alpha ,\ {\left\Vert{\varphi }\right\Vert}_{L^{2}(\Omega
 ,\eta )}\leq C{\left\Vert{\alpha }\right\Vert}_{L^{2}(\Omega ,\mu )}.$ \ \par 
Hence	  $\displaystyle {\mathcal{L}}(\alpha )=(-1)^{p+q-1}{\left\langle{\omega
 ,\varphi }\right\rangle}$  and by Cauchy Schwarz inequality\ \par 
\quad \quad \quad \quad \quad 	 $\ \left\vert{{\mathcal{L}}(\alpha )}\right\vert \leq {\left\Vert{\omega
 }\right\Vert}_{L^{2}(\Omega ,\eta ^{-1})}{\left\Vert{\varphi
 }\right\Vert}_{L^{2}(\Omega ,\eta )}\leq C{\left\Vert{\omega
 }\right\Vert}_{L^{2}(\Omega ,\eta ^{-1})}{\left\Vert{\alpha
 }\right\Vert}_{L^{2}(\Omega ,\mu )},$ \ \par 
hence\ \par 
\quad \quad \quad \quad \quad 	 $\ \left\vert{{\mathcal{L}}(\alpha )}\right\vert \leq C{\left\Vert{\omega
 }\right\Vert}_{L^{2}(\Omega ,\eta ^{-1})}{\left\Vert{\alpha
 }\right\Vert}_{L^{2}(\Omega ,\ \mu )}.$ \ \par 
\quad \quad  	So we have that the norm of  ${\mathcal{L}}$  is bounded on
 the subspace of  $\bar \partial $  closed forms in  $L^{2}(\Omega
 ,\mu )$  by  $\ C{\left\Vert{\omega }\right\Vert}_{L^{2}(\Omega
 ,\eta ^{-1})}.$ \ \par 
\quad \quad  	We apply the Hahn-Banach theorem to extend  ${\mathcal{L}}$
  with the {\sl same} norm to {\sl all}   $(n-p,n-q+1)$  forms
 in  $L^{2}(\Omega ,\ \mu ).$  As in Serre's duality theorem
 (~\cite{Serre55}, p. 20) this is one of the main ingredient
 in the proof.\ \par 
\quad \quad  	This means, by the definition of currents, that there is a
  $(p,q-1)$  current  $u$  which represents the extended form
  ${\mathcal{L}},$  i.e.\ \par 
\quad \quad \quad \quad \quad 	 $\displaystyle \forall \alpha \in L_{(n-p,n-q+1)}^{2}(\Omega
 ,\mu ),\ \exists \varphi \in L^{2}(\Omega ,\eta )::{\left\langle{u,\alpha
 }\right\rangle}={\mathcal{L}}(\alpha )=(-1)^{p+q-1}{\left\langle{\omega
 ,\varphi }\right\rangle},$ \ \par 
and such that\ \par 
\quad \quad \quad \quad \quad 	 $\displaystyle \sup _{\alpha \in L^{2}(\Omega ,\mu ),\ {\left\Vert{\alpha
 }\right\Vert}_{L^{2}(\Omega ,\mu )}\leq 1}\ \left\vert{{\left\langle{u,\alpha
 }\right\rangle}}\right\vert \leq C{\left\Vert{\omega }\right\Vert}_{L^{2}(\Omega
 ,\eta ^{-1})}.$ \ \par 
 $\hfill\blacksquare $ \ \par 
\quad \quad  	Proof of the theorem~\ref{LpPoidsStein11} and theorem~\ref{LpPoidsStein13}.\
 \par 
Let  $\varphi $  be a strictly plurisubharmonic function in
 the Stein manifold  $\displaystyle \Omega .$  Let  $\displaystyle
 c_{\varphi }(z)$  be the smallest eigenvalue of  $\partial \bar
 \partial \varphi (z),$  then  $\displaystyle \forall z\in \Omega
 ,\ c_{\varphi }(z)>0.$  \ \par 
\quad \quad  	If  $\displaystyle \Omega $  is a pseudo convex domain in 
 ${\mathbb{C}}^{n},$  by lemma 4.4.1. of H\"ormander~\cite{Hormander73},
 p. 92, we have that, with  $\displaystyle \eta =e^{-\varphi
 },\ \mu =c_{\varphi }^{-1}e^{-\varphi },$ \ \par 
\quad \quad \quad \quad \quad 	 $\forall \alpha \in L_{(n-p,n-q+1)}^{2}(\Omega ,\mu ),\ \bar
 \partial \alpha =0,\ \exists \varphi \in L_{(n-p,n-q)}^{2}(\Omega
 ,\eta )::\bar \partial \varphi =\alpha ,$ \ \par 
and by theorem 4.2.2. still from H\"ormander~\cite{Hormander73},
 p. 84, because if  $\displaystyle \beta \in L_{(n-p,n-q)}^{2}(\Omega
 ,\eta )$  then  $\displaystyle \beta \in L_{(n-p,n-q)}^{2}(\Omega
 ,loc),$  we have\ \par 
\quad \quad \quad \quad \quad 	 $\displaystyle \forall \beta \in L_{(n-p,n-q)}^{2}(\Omega ,\eta
 ),\ \bar \partial \beta =0,\ \exists \gamma \in L_{(n-p,n-q)}^{2}(\Omega
 ,loc)::\bar \partial \gamma =\beta .$ \ \par 
Hence 	it remains to apply lemma~\ref{LpPoidsStein10}, proved
 in the next section, with  $\displaystyle r=2$  and  $\eta =c_{\varphi
 }^{-1}e^{-\varphi }$  to get the theorem~\ref{LpPoidsStein11}.\ \par 
\quad \quad  	If  $\displaystyle \Omega $  is a Stein manifold, we need 
 $\displaystyle c_{\varphi }$  big enough to verify inequality
 (5.2.13), p. 125 in~\cite{Hormander73}. To have this we can
 replace  $\varphi $  by  $\psi :=\chi (\varphi )$  with  $\chi
 $  a convex increasing function, as in~\cite{Hormander73}, p.
 125. Then we still can apply what precedes to have the theorem~\ref{LpPoidsStein13}.
  $\hfill\blacksquare $ \ \par 

\subsection{A duality lemma with weight.}
\ \par 
\quad \quad  	This is still in~\cite{AnGrauLr12} lemma 3.1., but I repeat
 it completely for the reader's convenience.\ \par 
\quad \quad  	Let  ${\mathcal{I}}_{p}$  be the set of multi-indices of length
  $p$  in  $(1,...,n).$  	We shall use the measure defined on
  $\Gamma :=\Omega {\times}{\mathcal{I}}_{p}{\times}{\mathcal{I}}_{q}$
  the following way :\ \par 
\quad \quad \quad \quad   $\displaystyle d\mu (z,k,l)=d\mu _{\eta ,p,q}(z,k,l):=\eta
 (z)dm(z)\otimes \ \sum_{\left\vert{I}\right\vert =p,\ \left\vert{J}\right\vert
 =q}{\delta _{I}(k)\otimes \delta _{J}(l)},$ \ \par 
where  $\delta _{I}(k)=1$  if the multi-index  $k$  is equal
 to  $I$  and  $\delta _{I}(k)=0$  if not.\ \par 
This means, if  $f(z,I,J)$  is a function defined on  $\Gamma ,$  that\ \par 
\quad \quad \quad \quad \quad 	 $\displaystyle \ \int{f(z,k,l)d\mu _{\eta ,p,q}(z,k,l)}:=\sum_{\left\vert{I}\right\vert
 =p,\ \left\vert{J}\right\vert =q}{\int_{\Omega }{f(z,I,J)\eta
 (z)dm(z)}}.$ \ \par 
\quad \quad  	If  $I$  is a multi-index of length  $p,$  let  $I^{c}$  be
 the unique multi-index, ordered increasingly, such that  $I\cup
 I^{c}=(1,\ 2,...,\ n)\ ;$  then  $I^{c}$  is of length  $n-p.$ \ \par 
\quad \quad  	To  $t=\ \sum_{\left\vert{I}\right\vert =p,\ \left\vert{J}\right\vert
 =q}{t_{I,J}(z)dz^{I}\wedge d\bar z^{J}}$  a  $(p,q)$  form,
 we associate the function on  $\Gamma $  :\ \par 
\quad \quad \quad \quad \quad 	 $\displaystyle T(z,I,J):=(-1)^{s(I,J)}t_{I,J}(z),$ \ \par 
where\ \par 
\quad \quad \quad \quad \quad   $s(I,J)=0$  if  $dz^{I}\wedge d\bar z^{J}\wedge dz^{I^{c}}\wedge
 d\bar z^{J^{c}}=dz_{1}\wedge \cdot \cdot \cdot \wedge dz_{n}\wedge
 d\bar z_{1}\wedge \cdot \cdot \cdot \wedge d\bar z_{n}$  as
 a  $(n,n)$  form\ \par 
and\ \par 
\quad \quad \quad \quad \quad 	 $s(I,J)=1$  if not.\ \par 
\quad \quad  	If  $\varphi =\ \sum_{\left\vert{I}\right\vert =p,\ \left\vert{J}\right\vert
 =q}{\varphi _{I^{c},J^{c}}(z)dz^{I}\wedge d\bar z^{J}}$  is
 of complementary bi-degree, associate in the same manner :\ \par 
\quad \quad \quad \quad \quad 	 $\displaystyle \Phi ^{*}(z,\ I,\ J):=\varphi _{I^{c},J^{c}}(z).$
  This is still a function on  $\Gamma .$ \ \par 
\quad \quad  	Now we have, for  $1<r<\infty ,$  if  $T(z,I,J)$  is a function
 in  $\Omega $  with  $L^{r}(\Omega )$  coefficients and with
  $\mu =\mu _{\eta ,p,q},$ \ \par 
\quad \quad \quad \quad \quad 	 $\displaystyle \ {\left\Vert{T}\right\Vert}_{L^{r}(d\mu )}^{r}:=\int{\left\vert{T(z,I,J)}\right\vert
 ^{r}d\mu _{\eta ,p,q}(x,\ I,\ J)}=\ \ \sum_{\left\vert{I}\right\vert
 =p,\ \left\vert{J}\right\vert =q}{{\left\Vert{T(\cdot ,I,J)}\right\Vert}_{L^{r}(\Omega
 ,\eta )}^{r}}.$ \ \par 
\quad \quad  	For  $1\leq r<\infty $  the dual of  $L^{r}(\mu )$  is  $L^{r'}(\mu
 )$  where  $r'$  is the conjugate of  $\displaystyle r,\ \frac{1}{r}+\frac{1}{r'}=1,$
  and the norm is defined analogously with  $r'$  replacing  $r.$ \ \par 
\quad \quad  	We also know that, for  $p,\ q$  fixed,\ \par 
\begin{equation}  \ {\left\Vert{T}\right\Vert}_{L^{r}(\mu )}=\sup
 _{\Phi \in L^{r'}(\mu )}\frac{\left\vert{\int{T\Phi d\mu }}\right\vert
 }{{\left\Vert{\Phi }\right\Vert}_{L^{r'}(\mu )}}.\label{AndGrau127}\end{equation}\
 \par 
\quad \quad  	For a  $(p,\ q)$  form			  $t=\sum_{\left\vert{J}\right\vert
 =p,\ \left\vert{K}\right\vert =q}{t_{J,K}dz^{J}\wedge d\bar
 z^{K}},$  and a weight  $\eta >0$  we define its norm by :\ \par 
\begin{equation}  \ {\left\Vert{t}\right\Vert}_{L^{r}(\Omega
 ,\eta )}^{r}:=\ \ \sum_{\left\vert{J}\right\vert =p,\ \left\vert{K}\right\vert
 =q}{{\left\Vert{t_{J,K}}\right\Vert}_{L^{r}(\Omega ,\eta )}^{r}}={\left\Vert{T}\right\Vert}_{L^{r}(\mu
 )}^{r}.\label{AndGrau137}\end{equation}\ \par 
\quad \quad  	Now we can state\ \par 
\begin{Lemma}
~\label{LpPoidsStein10}Let  $\eta >0$  be a weight. If  $u$
  is a  $(p,q)$  current defined on  $(n-p,n-q)$  forms in  $L^{r'}(\Omega
 ,\eta )$  and such that\par 
\quad \quad \quad \quad \quad 	 $\forall \alpha \in L^{r'}_{(n-p,n-q)}(\Omega ,\eta ),\ \left\vert{{\left\langle{u,\
 \alpha }\right\rangle}}\right\vert \leq C{\left\Vert{\alpha
 }\right\Vert}_{L^{r'}(\Omega ,\eta )},$ \par 
then  $\ {\left\Vert{u}\right\Vert}_{L^{r}(\Omega ,\eta ^{1-r})}\leq C.$ \par 
\end{Lemma}
\quad \quad  	Proof.\ \par 
	Let us take the measure  $\mu =\mu _{\eta ,p,q}$  as above.
 Let  $\Phi ^{*}$  be the function on  $\Gamma $  associated
 to  $\alpha $  and  $T$  the one associated to  $u.$  We have,
 by definition of the measure  $\mu $  applied to the function\ \par 
\quad \quad \quad \quad \quad 	 $f(z,I,J):=T(z,I,J)\eta ^{-1}\Phi ^{*}(z,I,J),$  \ \par 
\quad \quad \quad \quad \quad 	 $\displaystyle \ \int{T\eta ^{-1}\Phi ^{*}d\mu }=\int{f(z,\
 k,\ l)d\mu (z,k,l)}:=\sum_{\left\vert{I}\right\vert =p,\ \left\vert{J}\right\vert
 =q}{\int_{\Omega }{f(z,\ I,\ J)\eta (z)dm(z)}}=$ \ \par 
\quad \quad \quad \quad \quad \quad \quad \quad \quad \quad \quad \quad \quad 	 $\displaystyle =\sum_{\left\vert{I}\right\vert =p,\ \left\vert{J}\right\vert
 =q}{\int_{\Omega }{T(z,I,J)\eta ^{-1}(z)\Phi ^{*}(z,I,J)\eta
 (z)dm(z)}}={\left\langle{u,\ \alpha }\right\rangle},$ \ \par 
by definition of  $T$  and  $\Phi ^{*}.$ \ \par 
Hence we have, by~(\ref{AndGrau127})\ \par 
\quad \quad \quad \quad \quad 	 $\displaystyle \ {\left\Vert{T\eta ^{-1}}\right\Vert}_{L^{r}(\mu
 )}=\sup _{\Psi \in L^{r'}(\mu )}\frac{\left\vert{{\left\langle{u,\
 \alpha }\right\rangle}}\right\vert }{{\left\Vert{\Psi }\right\Vert}_{L^{r'}(\mu
 )}}.$ \ \par 
But  $\ {\left\Vert{T\eta ^{-1}}\right\Vert}_{L^{r}(\mu )}={\left\Vert{u\eta
 ^{-1}}\right\Vert}_{L^{r}(\Omega ,\ \eta )}$  by~(\ref{AndGrau137}), and\ \par 
\quad \quad \quad \quad \quad 	 $\displaystyle \ {\left\Vert{f\eta ^{-1}}\right\Vert}_{L^{r}(\Omega
 ,\eta )}^{r}=\int_{\Omega }{\left\vert{f\eta ^{-1}}\right\vert
 ^{r}\eta dm}=\int_{\Omega }{\left\vert{f}\right\vert ^{r}\eta
 ^{1-r}dm}={\left\Vert{f}\right\Vert}_{L^{r}(\Omega ,\eta ^{1-r})},$ \ \par 
so we get\ \par 
\quad \quad \quad \quad \quad 	 $\displaystyle \ {\left\Vert{u}\right\Vert}_{L^{r}(\Omega ,\eta
 ^{1-r})}=\sup _{\Psi \in L^{r'}(\mu )}\frac{\left\vert{{\left\langle{u,\
 \alpha }\right\rangle}}\right\vert }{{\left\Vert{\Psi }\right\Vert}_{L^{r'}(\mu
 )}},$ \ \par 
which implies the lemma because, still by~(\ref{AndGrau127}),
 we can take  $\Psi =\Phi ^{*}$  and  $\ {\left\Vert{\Psi }\right\Vert}_{L^{r'}(\mu
 )}={\left\Vert{\alpha }\right\Vert}_{L^{r'}(\Omega ,\eta )}.$
   $\hfill\blacksquare $ \ \par 
\quad \quad  	It may seem strange that we have such an estimate when the
 dual of  $L^{r'}(\Omega ,\eta )$  is   $\displaystyle L^{r}(\Omega
 ,\eta ),$  but the reason is, of course, that in the duality
 forms-currents there is no weights.\ \par 
\ \par 

\bibliographystyle{C:/texlive/2012/texmf-dist/bibtex/bst/base/plain}

\end{document}